\documentclass[11pt]{article}%
\usepackage{amsfonts}
\usepackage{graphicx}
\usepackage{amsmath}
\usepackage{amssymb}
\usepackage{dsfont}
\usepackage[a4paper,bottom=3.5cm,top=2cm,left=2.5cm,right=2.5cm]{geometry}%
\providecommand{\U}[1]{\protect\rule{.1in}{.1in}}

\newtheorem{theo}{Theorem}[section]
\newtheorem{prop}[theo]{Proposition}

\newtheorem{coro}[theo]{Corollary}

\newtheorem{Rq}[theo]{Remark}

\newcommand\beq{\begin{equation}}
\newcommand\eeq{\end{equation}}
\newcommand{\E}{{\ensuremath{\mathbb E}} }

\def\nbX0{\widehat{X}_{0,m}}

\def\ntX0{\widetilde{X}_{0,m}}

\begin{document}
\begin{center} {\bf \Large Limiting spectral distribution of Gram matrices associated with functionals of $\beta$-mixing processes }\vskip15pt

Marwa Banna 
\let\thefootnote\relax\footnote{Universit\'e Paris Est, LAMA (UMR 8050), UPEMLV, CNRS, UPEC, 5 Boulevard Descartes, 77454 Marne La
Vall\'ee, France. \\ 
E-mail: marwa.banna@u-pem.fr\vskip10pt}
\end{center}

\begin{abstract} We give asymptotic spectral results for Gram matrices of the form $ n^{-1}\mathcal{X}_n \mathcal{X}_n^T$ where the entries of $\mathcal{X}_n$ are dependent across both rows and columns. More precisely, they consist of short or long range dependent random variables having moments of second order and that are functionals of an absolutely regular sequence. We also give a concentration inequality of the Stieltjes transform and we prove that, under an arithmetical decay condition on the $\beta$-mixing coefficients, it is almost surely concentrated around its expectation. Applications to examples of positive recurrent Markov chains and dynamical systems are also given.
\end{abstract}



\section{Introduction}

\quad
For a random matrix $\mathcal{X}_n \in \mathbb{R}^{N \times n}$, the study of the asymptotic behavior of the eigenvalues of the $N \times N$ Gram matrix $n^{-1}\mathcal{X}_n \mathcal{X}_n^T$ gained interest as it is employed in many applications in statistics, signal processing, quantum physics, finance, etc. In order to describe the distribution of the eigenvalues, it is convenient to introduce the empirical spectral measure defined by
$
\mu_{n^{-1}\mathcal{X}_n \mathcal{X}_n^T}= N^{-1} \sum_{k=1}^{N} \delta_{\lambda_k} \, ,
$
where $\lambda_1, \ldots , \lambda_N$ are the eigenvalues of $n^{-1}\mathcal{X}_n \mathcal{X}_n^T$. This type of study was actively developed after the pioneering work of Mar\u{c}enko and Pastur \cite{MaPa}, who proved that under the assumption $\lim_{n \rightarrow + \infty} N/n = c \in (0,+\infty)$, the empirical spectral distribution of large dimensional Gram matrices with i.i.d. centered entries having finite variance converges almost surely to a non-random distribution. The limiting spectral distribution (LSD) obtained, i.e. the Mar\u{c}enko-Pastur distribution, is given explicitly in terms of $c$ and depends on the distribution of the entries of $\mathcal{X}_n$ only through their common variance. The original Mar\u{c}enko-Pastur theorem is stated for random variables having moment of fourth order; for the proof under second moment only, we refer to Yin \cite{Yi}. 

\smallskip
Since then, a large amount of study has been done aiming to relax the independence structure between the entries of $\mathcal{X}_n$. For example, Bai and Zhou \cite{BZ} treated the case where the columns of $\mathcal{X}_n$ are independent with their coordinates having a very general dependence structure and moments of fourth order. Recently, Banna and Merlev\`ede \cite{BaMe} extended along another direction the Mar\u{c}enko-Pastur theorem to a large class of weakly dependent sequences of real random variables having moments of second order. Letting $(X_k)_{k \in {\mathbb Z}}$ be a stationary process of the form $X_k=g(\cdots, \varepsilon_{k-1}, \varepsilon_k, \varepsilon_{k+1}, \ldots )$, where the $\varepsilon_k$'s are i.i.d. real valued random variables and $g:{\mathbb R}^{{\mathbb Z}} \rightarrow {\mathbb R}$ is a measurable function, they consider the $N \times N$ sample covariance matrix ${\bf A}_n = \frac{1}{n} \sum_{k=1}^{n} {\bf X}_k {\bf X}_k^T$ with the ${\bf X}_k$'s being independent copies of the vector ${\bf X}= (X_1, \dots ,X_N)^T$. Assuming only that $X_0$ has a moment of second order, then provided that $\lim_{n \rightarrow \infty} N/n = c \in (0,\infty)$, they prove, under a mild dependence condition, that almost surely, $\mu_{{\bf A}_n}$ converges weakly to a non-random probability measure $\mu$ whose Stieltjes transform satisfies an integral equation that depends on $c$ and on the spectral density of the underlying stationary process $(X_k)_{k \in {\mathbb Z}}$. In this paper, we relax the independence condition on the $\varepsilon_k$'s. For instance, assuming that $(\varepsilon_k)_{k \in \mathbb{Z}}$ is an absolutely regular ($\beta$-mixing) sequence, we show that the study of the LSD can be reduced to that of a Gaussian matrix having the same covariance structure (see Theorem \ref{second theo}). Then, provided that the spectral density of the underlying process exists, we deduce the integral equation satisfied by the Stieltjes transform of the LSD $\mu$.

\smallskip
In the above mentioned model, the random vector ${\bf X}= (X_1, \dots ,X_N)^T$ can be viewed as an $N$-dimensional process repeated independently $n$ times to obtain the ${\bf X}_k$'s. However, in practice it is not always possible to observe a high dimensional process several times. In the case where only one observation of length $Nn$ can be recorded, it seems reasonable to partition it into $n$ dependent observations of length $N$, and to treat them as $n$ dependent observations. Up to our knowledge this was first done by Pfaffel and Schlemm \cite{PfSc} who showed that this approach is valid and leads to the correct asymptotic eigenvalue distribution of the sample covariance matrix if the components of the underlying process are modeled as short memory linear filters of independent random variables. They consider Gram matrices having the same form as in \eqref{matrix-Bn} and associated with a stationary linear process $(X_k)_{k \in \mathbb{Z}}$ with independent innovations having finite fourth moments and such that the coefficients decay with an arithmetical rate, and they derive its LSD.

\smallskip
In this work, we study the same model of random matrices as in \cite{PfSc} but considering the case where the entries come from a non causal stationary process $(X_k)_{k \in \mathbb{Z}}$ of the form $X_k=g(\cdots, \varepsilon_{k-1}, \varepsilon_k , \varepsilon_{k+1}, \ldots )$ where $(\varepsilon_k)_{k \in \mathbb{Z}}$ is an absolutely regular sequence and $g:{\mathbb R}^{{\mathbb Z}} \rightarrow {\mathbb R}$ is a measurable function such that $X_k$ is a proper centered random variable having finite moments of second order. Under an arithmetical decay condition on the $\beta$-mixing coefficients, we prove in Theorem \ref{theo conc ineq} that the Stieltjes transform is concentrated almost surely around its expectation as $n$ tends to infinity (see Theorem \ref{theo conc ineq}). Once this is done, it is enough to prove that the expectation of the Stieltjes transform converges to that of a non-random probability measure.  This can be achieved by approximating it by the expectation of the Stieltjes transform of a Gaussian matrix having a close covariance structure (see Theorem \ref{first theorem}). Finally, provided that the spectral density of $(X_k)_k$ exists, we prove in Theorem \ref{full theo} that almost surely, $\mu_{{\bf B}_n}$ converges weakly to the same non-random limiting probability measure $\mu$ obtained in the cases mentioned above.

We recall now that the absolutely regular ($\beta$-mixing) coefficient between two $\sigma$-algebras $\mathcal{A}$ and $\mathcal{B}$ is defined by
$$
\beta ( \mathcal{A} ,  \mathcal{B})=  \frac{1}{2} \sup \Big\lbrace \sum_{i \in I} \sum_{j \in J} \big|\mathbb{P} (A_i \cap B_j) - \mathbb{P}(A_i) \mathbb{P}(B_j) \big| \Big\rbrace \, ,
$$
where the supremum is taken over all finite partitions $(A_i)_{i \in I}$ and $(B_j)_{j \in J}$ that are respectively $\mathcal{A}$ and $\mathcal{B}$ measurable (see Rozanov and Volkonskii \cite{VR}). The coefficients $(\beta_n)_{n \geqslant 0}$ of a sequence $(\varepsilon_i)_{i \in \mathbb{Z}}$ are defined by 
\begin{equation}\label{def beta n}
\beta_0 = 1 \quad \text{and} \quad \beta_n = \sup_{k \in \mathbb{Z}} \beta \big( \sigma (\varepsilon_{\ell} \: , \: \ell \leqslant k) \; , \; (\varepsilon_{\ell +n} \: , \: \ell \geqslant k) \big) \quad \text{for} \quad n\geqslant 1.
\end{equation}
Moreover, $(\varepsilon_i)_{i \in \mathbb{Z}}$ is said to be absolutely regular or $\beta$-mixing if $\beta_n \rightarrow 0$ as $n \rightarrow \infty$.

\smallskip
\noindent
\textbf{Outline.} In Section 2, we specify the models studied and state the limiting results for the Gram matrices associated with the process defined in \eqref{def Xn}. The proofs shall be deferred to Section 4, whereas applications to examples of Markov chains and dynamical systems shall be introduced in Section 3.

\smallskip
\noindent
\textbf{Notation.}
For any real numbers $x$ and $y$, $ x \wedge y:=\min (x,y)$ whereas $x \vee y:=\max (x,y)$. Moreover, the notation $[x]$ denotes the integer part of $x$. For any non-negative integer $q$, a null row vector of dimension $q$ will be denoted by ${\bf 0}_{q}$. For a matrix $A$, we denote by $A^T$ its transpose matrix and by ${\rm Tr} (A)$ its trace.  Finally, we shall use the notation $\Vert X \Vert_r$ for the ${\mathbb L}^r$-norm ($r \geq 1$) of a real valued random variable $X$.

For any square matrix $A$ of order $N$ having only real eigenvalues, its empirical spectral measure and distribution are respectively defined by
$$
\mu_A=\frac{1}{N} \sum_{k=1}^{N} \delta_{\lambda_k}
\ {\rm and} \
F^{A}(x) =\frac{1}{N} \sum_{k=1}^{N}  {\bf 1}_{ \{ \lambda_k \leqslant x \}} \, , 
$$
where $\lambda_1, \ldots , \lambda_N$ are the eigenvalues of $A$. The Stieltjes transform of $\mu_{A}$ is given by 
$$
 S_A (z) := S_{\mu_{A}}(z)= \int \frac{1}{x-z} d\mu_{A}(x)= \frac{1}{N} {\rm Tr}(A - z{\bf I})^{-1} \, ,
$$
where $z=u+iv\in \mathbb{C}^+$ (the set of complex numbers with positive imaginary part), and $\bf{I}$ is the identity matrix of order $N$.

\section{Results}\label{section results}
\setcounter{equation}{0}

\quad We consider a non causal stationary process $(X_k)_{k \in \mathbb{Z}}$ defined as follows: let  $(\varepsilon_i)_{i \in \mathbb{Z}}$ an absolutely regular process with $\beta$-mixing sequence $(\beta_k)_{k\geqslant 0}$ and let $g: \mathbb{R}^{\mathbb{Z}} \rightarrow \mathbb{R} $ be a measurable function such that $X_k$, which is defined for any $k \in \mathbb{Z}$ by
\begin{equation}\label{def Xn}
X_k=g (\xi_k) \quad \text{with} \quad \xi_k= (\ldots, \varepsilon_{k-1}, \varepsilon_k, \varepsilon_{k+1} , \ldots),
\end{equation}
is a proper centered random variable having finite moment of second order; that is, $\E (X_k)=0$ and $\Vert X_k \Vert_2 <\infty$.

\noindent
Now, let $N:=N(n)$ be a sequence of positive integers and consider the $N \times n$ random matrix $\mathcal{X}_n$ defined by
\begin{equation}\label{matrix-Xn}
\mathcal{X}_n=((\mathcal{X}_n)_{i,j})= (X_{(j-1)N+i})=
\left(
\begin{array}{cccc}
X_1     &  X_{N+1}  & \cdots  & X_{(n-1)N+1} \\
X_2     &  X_{N+2}  & \cdots  & X_{(n-1)N+2} \\
\vdots  &  \vdots   & \,	   & \vdots	 \\
X_N     &  X_{2N}   & \cdots  & X_{nN}
\end{array}
 \right)
\in \mathcal{M}_{N \times n} (\mathbb{R})
 \end{equation}
and note that its entries are dependent across both rows and columns. Let ${\bf B}_n$ be its corresponding Gram matrix given by
\begin{equation} \label{matrix-Bn}
{\bf B}_n= \frac{1}{n} \mathcal{X}_n \mathcal{X}_n^T.
\end{equation} 
In what follows, ${\bf B}_n$ will be referred to as the Gram matrix associated with $(X_k)_{k\in \mathbb{Z}}$. We note that it can be written as the sum of dependent rank one matrices. Namely, ${\bf B}_n=\frac{1}{n} \sum_{i=1}^{n}{\bf X}_i{\bf X}_i ^T$, where for any $i=1,\ldots, n$, ${\bf X}_i= (X_{(i-1)N+1}, \ldots, X_{iN})^T$. Our purpose is to study the limiting distribution of the empirical spectral measure $\mu_{{\bf B}_n}$ defined by
$$
\mu_{{\bf B}_n}(x)= \frac{1}{N} \sum_{k=1}^{N} \delta_{\lambda_k} \, ,
$$
where $\lambda_1, \ldots , \lambda_N$ are the eigenvalues of ${\bf B}_n$. We start by showing that if the $\beta$-mixing coefficients decay arithmetically then the Stieltjes transform of ${\bf B}_n$ concentrates almost surely around its expectation as $n$ tends to infinity .

\begin{theo}\label{theo conc ineq}
Let ${\bf B}_n$ be the matrix defined in \eqref{matrix-Bn} and associated with $(X_k)_{k \in \mathbb{Z}}$ defined in \eqref{def Xn}. If $\lim_{n\rightarrow \infty} N/n = c \in (0, \infty)$ and
\beq \label{Cond-beta}
\sum_{n \geqslant 1 } \log (n)^{\frac{3 \alpha}{2}} n^{-\frac{1}{2}} \beta_n < \infty  \ \ \text{for some} \ \alpha >1 \, ,
\eeq
the following convergence holds: for any $z \in \mathbb{C}^{+}$,
$$
 S_{{{\bf B}_n}}(z) - \E  \big( S_{{{\bf B}_n}}(z)\big)\rightarrow 0 \ \text{almost surely, as} \ n \rightarrow +\infty \, . 
$$
\end{theo}
As the ${\bf X}_i$'s are dependent, classical arguments as those in Theorem 1 $(ii)$ of \cite{GuLe} are not sufficient to prove the above convergence. In fact, we use maximal coupling for absolutely regular sequences in order to break the dependence structure between the columns. The proof of Theorem \ref{theo conc ineq} shall be postponed to Section \ref{section-1}.

\medskip

In the following Theorem, we shall approximate the Stieltjes transform of the LSD of ${\bf B}_n$ with that of a sample covariance matrix ${\bf G}_n$ which is the sum of i.i.d. rank one matrices associated with a Gaussian process $(Z_k)_{ k\in \mathbb{Z}}$ having the same covariance structure as $(X_k)_{ k\in \mathbb{Z}}$. Namely, for any $k,\ell \in \mathbb{Z}$,
\beq \label{covG} 
{\rm Cov}(Z_k,Z_{\ell})={\rm Cov}(X_k, X_{\ell}) \, .
\eeq
Denoting, for  $i = 1, \dots, n$, by $( Z^{(i)}_k)_{ k\in \mathbb{Z}}$ an independent copy of $(Z_k)_{k}$ that is also independent of $(X_k)_{k}$, we  then define  the $N \times N$ sample covariance matrix ${\bf G}_n$ by
\begin{equation}\label{def Gn}
{\bf G}_n=\frac{1}{n}\mathcal{Z}_n\mathcal{Z}_{n}^{T} \,= \frac{1}{n} \sum_{k=1}^{n} {\bf Z}_i {\bf Z}_i^T \,,
\end{equation} 
where for any $i=1, \ldots ,n$, ${\bf Z}_i=(Z_1^{(i)}, \ldots, Z_N^{(i)})^T$ and ${\mathcal{Z}}_n$ is the matrix whose columns are the ${\bf Z}_i$'s. Namely, ${\mathcal{Z}}_n:= (({\mathcal{Z}}_n)_{u,v})=(Z_u^{(v)})$. 

\begin{theo}\label{first theorem}
Let ${\bf B}_n$ and ${\bf G}_n$ be the matrices defined in \eqref{matrix-Bn} and \eqref{def Gn} respectively. Provided that $\lim_{n\rightarrow \infty} N/n = c \in (0, \infty)$ and $\lim_{n\rightarrow \infty}\beta_n=0 $, then for any $z \in \mathbb{C}_{+}$,
$$
\lim_{n\rightarrow \infty} | \E \big(S_{{\bf B }_n}(z) \big) - \E \big(S_{{\bf G}_n}(z)\big) | =0. 
$$ 
\end{theo}
The above Theorem allows us to reduce the study of the expectation of the LSD of $\textbf{B}_n$ to that of a Gram matrix, being the sum of independent rank-one matrices associated with a Gaussian process, without requiring any rate of convergence to zero of the correlation between the entries nor of the $\beta$-mixing coefficients.

\begin{theo}\label{full theo}
Let ${\bf B}_n$ and ${\bf G}_n$ be the matrices defined in \eqref{matrix-Bn} and \eqref{def Gn} respectively. Provided that $\lim_{n\rightarrow \infty} N/n = c \in (0, \infty)$ and that \eqref{Cond-beta} is satisfied, then for any $z \in \mathbb{C}_{+}$,
\begin{equation}\label{Conv-main-theo}
\lim_{n\rightarrow \infty} | S_{{\bf B }_n}(z)  -S_{{\bf G}_n}(z) | =0 \quad \text{a.s.} 
\end{equation}
Moreover, if $(X_k)_{ k\in \mathbb{Z}}$ admits a spectral density $f$, then with probability one, $\mu_{{\bf B}_n}$ converges weakly to a probability measure $\mu$ whose Stieltjes transform $S=S(z)$ ($z \in \mathbb{C}^+$) satisfies the equation 
\beq  \label{equ S(z)} 
z = -\frac{1}{{\underline S}} + \frac{c}{2  \pi} \int_0^{2\pi} \frac{1}{{\underline S} + \big ( 2 \pi f(\lambda) \big )^{-1}} d \lambda \, ,
\eeq
where ${\underline S}(z):=-(1-c)/z +c S(z)$.
\end{theo}

\begin{Rq} \label{Rq-spectral-density} The spectral density function $f$ of $(X_k)_{k \in {\mathbb Z}}$ is  the discrete Fourier transform of the autocovariance function. If $\sum_{k \in \mathbb{Z}} |\text{Cov}(X_0, X_k)|<\infty$ then $f$ exists, is continuous and bounded on $[0,2\pi)$. It also follows from Proposition 1 in \cite{Ya} that the LSD $\mu$ is compactly supported. 
\end{Rq}

\begin{Rq}
For a proof of Equation \eqref{equ S(z)}, we refer to Section 4.4 in \cite{BaMe} for the case of short memory dependent sequences and to Section 3.5 in \cite{MePe} for the case of long memory dependent sequences adapted to the natural filtration. We also note that the convergence \eqref{Conv-main-theo} is a consequence of Theorems \ref{theo conc ineq} and \ref{first theorem} and therefore we omit the proof of Theorem \ref{full theo}.

\end{Rq}

Now, we shall consider the case where the Gram matrix is the sum of independent rank-one matrices whose entries are functionals of absolutely regular sequences. More precisely, we consider the sample covariance matrix 
\begin{equation}\label{def An}
{\bf A}_n:=\frac{1}{n}\sum_{k=1 }^{n} \widetilde{\bf X}_{k} \widetilde{\bf X}_{k}^{T},
\end{equation}
with the  $\widetilde{\bf X}_{k}$'s being independent copies of $(X_1, \ldots , X_N)^T$ and $(X_k)_{k \in \mathbb{Z}}$ the process defined in \eqref{def Xn}.
\begin{theo}\label{second theo}
Provided that $\lim_{n\rightarrow \infty} N/n = c \in (0, \infty)$ and $\lim_{n \rightarrow \infty}\beta_n=0 $, then for any $z \in \mathbb{C}$
\begin{equation}\label{conv-second-theo}
\lim_{n\rightarrow \infty} | S_{{\bf A }_n}(z)  -S_{{\bf G}_n}(z) | =0 \quad \text{a.s.} 
\end{equation}
Moreover, if $(X_k)_{ k\in \mathbb{Z}}$ admits a spectral density $f$,
then, with probability one, $\mu_{{\bf A}_n}$ converges weakly to a probability measure whose Stieltjes transform $S=S(z)$ ($z \in \mathbb{C}^+$) satisfies equation \eqref{equ S(z)}. 
\end{theo}

\begin{Rq}
Since the $\widetilde{\bf X}_{k}$'s are mutually independent, then, by Theorem 1 $(ii)$ in \cite{GuLe} or the arguments on page 34 in \cite{BS}, we can approximate directly $S_{{{\bf A}_n}}(z)$ by its expectation and there is no need of any coupling arguments as in Theorem \ref{theo conc ineq} and thus of the arithmetic decay condition $\eqref{Cond-beta}$ on the absolutely regular coefficients. So it suffices to prove that for any $z \in \mathbb{C}_{+}$,
$$
\lim_{n\rightarrow \infty} | \E \big(S_{{\bf A }_n}(z) \big) -\E \big(S_{{\bf G}_n}(z) \big)| =0 \,  
$$  
which can be exactly done as in Theorem \ref{first theorem} after simple modifications of indices. For this reason, its proof shall be omitted.
\end{Rq}

\section{Applications}

In this section we shall apply the results of Section \ref{section results} to a Harris recurrent Markov chain and some uniformly expanding maps in dynamical systems.

\subsection*{Harris recurrent Markov chain}
The following example is a symmetrized version of the Harris recurrent Markov chain defined by Doukhan \textit{et al.} \cite{DMR}. Let $(\varepsilon_n)_{n \in \mathbb{Z}}$ be a stationary Markov chain taking values in $E=[-1, 1]$ and let $K$ be its Markov kernel defined by
$$
K(x,.)=(1-|x|)\delta_x +|x| \nu \, ,
$$ 
with $\nu$ being a symmetric atomless law on $E$ and $\delta_x$ denoting the Dirac measure at point $x$. Assume that $ \theta = \int_E |x|^{-1} \nu (dx) < \infty$ then $(\varepsilon_n)_{n \in \mathbb{Z}}$ is positively recurrent and the unique invariant measure $\pi$ is given by
$$
\pi(dx)= \theta^{-1} |x|^{-1}\nu (dx)\, .
$$
We shall assume in what follows that $\nu$ satisfies for any $x \in [0,1]$,
\beq \label{assmptn nu}
 \frac{d \nu}{dx}(x) \leqslant c\, x^{a} \quad {\rm for \, some }\,\, a \,,\,c >0 \, .
 \eeq
Now, let $g$ be a measurable function defined on $E$ such that
\beq\label{defX app1}
X_k= g( \varepsilon_k)
\eeq
is a centered random variable having a finite second moment.
\begin{coro}
Let $(X_k)_{k \in \mathbb{Z} }$ be defined in \eqref{defX app1}. Assume that $\nu$ satisfies \eqref{assmptn nu} and that for any $x \in E$, $g(-x)=-g(x)$ and $|g(x)|\leqslant C |x|^{1/2}$ with $C$ being a positive constant.  Then, provided that $N/n \rightarrow c \in (0, \infty)$, the conclusions of Theorems \ref{first theorem} and \ref{second theo} hold and the LSD $\mu$ has a compact support. In addition, if \eqref{assmptn nu} holds with $a>1/2$ then Theorems \ref{theo conc ineq} and \ref{full theo} follow as well.
\end{coro}

\smallskip
\noindent
\textit{Proof.} Doukhan \textit{et al.} prove in Section 4 of \cite{DMR} that if \eqref{assmptn nu} is satisfied then $(\varepsilon_k)_{k \in \mathbb{Z}}$ is an absolutely regular sequence with $\beta_n= O(n^{-a})$ as $n\rightarrow \infty$. Thereby, Theorems \ref{first theorem} and \ref{second theo} follow. Now, noting that $g$ is an odd function we have
$$
\E (g(\varepsilon_k)|\varepsilon_0)= (1- |\varepsilon_0|)^k g(\varepsilon_0) \quad {\rm a.s.}
$$
Therefore, by the assumption on $g$ and \eqref{assmptn nu}, we get for any $k \geqslant 0$, 
\begin{align}
\gamma_k:=\E (X_0 X_k)
&=\E \big( g (\varepsilon_0) \E (g(\varepsilon_k)|\varepsilon_0)\big)
 = \theta^{-1} \int_E g^2(x) (1-|x|)^k |x|^{-1} \nu (dx) \notag \\
&\leqslant c\, C^2  \, \theta^{-1} \int_E \frac{x^{a+1}}{|x|} (1-|x|)^k  dx \, .
\end{align}
By the properties of the Beta and Gamma functions, $|\gamma_k|=O\big(\frac{1}{k^{a+1}}\big)$ which implies $\sum_k |\gamma_k|<\infty$ and thus the spectral density $f$ is continuous and bounded over $[0,2\pi)$ and the LSD $\mu$ has a compact support (see Remark \ref{Rq-spectral-density}).
However, if in addition $a>1/2$ then $\eqref{Cond-beta}$ is also satisfied and Theorems \ref{theo conc ineq} and \ref{full theo} follow as well.

\subsection*{Uniformly expanding maps}
Functionals of absolutely regular sequences occur naturally as orbits of chaotic dynamical systems. For instance, for uniformly expanding maps $T:[0,1] \rightarrow [0,1]$ with absolutely continuous invariant measure $\nu$, one can write $ T^k = g(\varepsilon_k, \varepsilon_{k+1} , \ldots)$ for some measurable function $g: \mathbb{R}^{\mathbb{Z}} \rightarrow \mathbb{R}$ where $(\varepsilon_k)_{k \geqslant 0}$ is an absolutely regular sequence. We refer to Section 2 of \cite{HoKe} for more details and for a precise definition of such maps (see also Example 1.4 in \cite{BoBuDe}). Hofbauer and Keller prove in Theorem 4 of \cite{HoKe} that the mixing rate of $(\varepsilon_k)_{k \geqslant 0}$ decreases exponentially, i.e. 
\beq\label{beta-m-c}
\beta_k \leqslant C e^{-\lambda k} \, , \quad {\rm for \ some \ } \, C , \, \lambda \,>0 \;, 
\eeq
and thus $\eqref{Cond-beta}$ holds. Setting for any $k \geqslant 0$,
\beq
X_k= h \circ T^k- \nu(f),
\eeq
where $h:[0,1] \rightarrow \mathbb{R}$ is a continuous H\"older function,  the Theorems in Section \ref{section results} hold for the associated matrices ${\bf B}_n$ and ${\bf A}_n$. Moreover,  Hofbauer and Keller prove in Theorem 5 of \cite{HoKe} that $\sum_k |\text{Cov}(X_0, X_k)|<\infty$ which implies that the spectral density $f$ exists, is continuous and bounded on $[0,2\pi)$ and that the LSD $\mu$ is compactly supported.

\section{Proof of Theorem \ref{theo conc ineq} } \label{section-1} 

\setcounter{equation}{0}

Let $ m $ be a positive integer (fixed for the moment) such that $m\leqslant \sqrt{N}/2$ and let $({X}_{k,m})_{k \in \mathbb{Z}}$ be the sequence defined for any $k \in \mathbb{Z}$ by, 
\beq
{X}_{k,m}= \E (X_k | \varepsilon_{k-m}, \ldots , \varepsilon_{k+m}).
\eeq
Consider the $N \times n $ matrix  $ {\mathcal{X}}_{n,m}=(({\mathcal{X}}_{n,m})_{i,j})=({X}_{(j-1)N+i, \, m})$ and finally set
\beq \label{matrix Bn,m bar}
{\bf B}_{n,m}=\frac{1}{n} {\mathcal{X}}_{n,m} {\mathcal{X}}_{n,m}^T.
\eeq
The proof will be done in two principal steps. First, we shall prove that 
\begin{equation}\label{app Bn by Bn,l}
\lim_{m \rightarrow \infty} \limsup_{n \rightarrow \infty}
\big| S_{{{\bf B}_n}}(z) - S_{{{\bf B}_{n,m}}}(z) \big| =0 \quad {\rm a.s.}
\end{equation}
and then
\begin{equation} \label{app Bn,l by its expc}
 \lim_{n \rightarrow \infty}
\big|  S_{{\bf B}_{n,m}}(z) - \E \big( S_{{\bf B}_{n,m}}(z) \big) \big| =0 \quad {\rm a.s.}
\end{equation}
We note that for any two $N \times n$ random matrices ${\bf A}$ and ${\bf B}$, we have
\beq\label{BaMe prop 4.1}
\big|S_{{{\bf A}{\bf A}^T}}(z) - S_{{{\bf B}{\bf B}^T}}(z) \big|
\leqslant \frac{\sqrt{2}}{N v^2} \left|  {\rm Tr}\big({\bf A}{\bf A}^T + {\bf B}{\bf B}^T \big) \right|^{1/2}
		  \left|  {\rm Tr} ( {\bf A} -{\bf B}) ( {\bf A} -{\bf B})^T \right|^{1/2}.
\eeq 
For a proof, the reader can check Inequalities (4.18) and (4.19) of \cite{BaMe}.
Thus, we get, for any $z=u+iv  \in \mathbb{C}^{+}$,
\begin{equation}\label{HW-WD 1}
\big|S_{{{\bf B}_n}}(z) - S_{{\bf B}_{n,m}}(z) \big|^2
\leqslant \frac{2}{v^4} \left( \frac{1}{N} {\rm Tr}({\bf B}_n + {\bf B}_{n,m} ) \right)
		  \left( \frac{1}{Nn} {\rm Tr} ( {\mathcal{X}}_n - \mathcal{X}_{n,m}) ( {\mathcal{X}}_n -\mathcal{X}_{n,m})^T \right)\, .
\end{equation}   
Recall that mixing implies ergodicity and note that as $(\varepsilon_k)_{k \in \mathbb{Z}}$ is an ergodic sequence of real-valued random variables then $(X_k)_{k \in \mathbb{Z}}$ is also so. Therefore, by the ergodic theorem,
\begin{equation}\label{ergth Xk}
\lim_{n \rightarrow + \infty} \frac{1}{N} {\rm Tr} ({\bf B}_n)
= \lim_{n \rightarrow + \infty} \frac{1}{Nn} \sum_{k =1}^{Nn}  X_{k}^2
= \E( X_0^2) \quad {\rm a.s.}
\end{equation}
Similarly,
\begin{align}
\lim_{n \rightarrow \infty} \frac{1}{N} {\rm Tr} ({\bf B}_{n,m})
= \E \big( X_{0, m}^2 \big) 
\quad {\rm a.s.} 
\end{align}
Starting from \eqref{HW-WD 1} and noticing that $\E( {X}_{0,m}^2) \leqslant \E (X_0^2)$, it follows that \eqref{app Bn by Bn,l} holds if we prove
\beq \label{CS app Bn by Bn bar}
\lim_{ m \rightarrow + \infty} \limsup_{n \rightarrow + \infty}  \left| \frac{1}{Nn} {\rm Tr} ( {\mathcal{X}}_n -\mathcal{X}_{n,m}) ( {\mathcal{X}}_n -{\mathcal{X}}_{n,m})^T \right| = 0 \quad {\rm a.s.}
\eeq
By the construction of $\mathcal{X}_n$ and ${\mathcal{X}}_{n,m}$ and again the ergodic theorem, we get
\begin{align*}
\lim_{n \rightarrow \infty}\frac{1}{Nn} {\rm Tr} (\mathcal{X}_n -{\mathcal{X}}_{n,m}) ( \mathcal{X}_n -{\mathcal{X}}_{n,m})^T
	 =\lim_{n \rightarrow \infty}  \frac{1}{Nn} \sum_{k =1}^{Nn}  (X_{k} - {X}_{k, \, m})^2
	= \E \big( X_{0} - {X}_{0, m}\big)^2 \quad {\rm a.s.} 
\end{align*}
\eqref{CS app Bn by Bn bar} follows by applying the usual martingale convergence theorem in $\mathbb{L}^2$, from which we infer that
$
\lim_{m \rightarrow + \infty} \Vert X_0 - \E( X_0| \varepsilon_{-m} , \ldots , \varepsilon_{m}) \Vert_2 = 0
$
(see Corollary 2.2 in Hall and Heyde \cite{HaHe}).

\smallskip
We turn now to the proof of \eqref{app Bn,l by its expc}. With this aim, we shall prove that for any $z=u+iv$ and $x>0$,
\beq \label{conc ineq stieltjes Bn,m bar}
\mathbb{P} \big( \big| S_{{\bf B}_{n,m}}(z) - \E S_{{\bf B}_{n,m}}(z)\big| > 4x\big) \leqslant 4 \exp\left\lbrace - \frac{x^2 \,v^2\, N^2 (\log n)^{\alpha}}{256  \, n^2 } \right\rbrace 
+ \frac{32 \, n^2 (\log n)^{\alpha}}{x^2 \,v^2\, N^2} \beta_{\big[\frac{n}{(\log n)^{\alpha}}\big]N}\, ,
\eeq
for some $\alpha >1$. Noting that
 $$
 \sum_{n \geqslant 2}\,(\log n)^{\alpha} \, \beta_{\big[\frac{n^2}{(\log n)^{\alpha}}\big]}  < +\infty \  \text{is equivalent to} \  \eqref{Cond-beta}
 $$ 
and applying Borel-Cantelli Lemma, \eqref{app Bn,l by its expc} follows  from \eqref{Cond-beta} and the fact that $ \lim_{n\rightarrow \infty} N/n= c$. 
Now, to prove \eqref{conc ineq stieltjes Bn,m bar}, we start by noting that 
\begin{multline}
\mathbb{P}\Big( \big|  S_{{{\bf B}_{n,m}}}(z) - \E \, \big(  S_{{{\bf B}_{n,m}}}(z)\big) \big| > 4x \Big) \notag \\
 \leqslant \mathbb{P}\Big( \!\big| \mathfrak{Re}\big( S_{{{\bf B}_{n,m}}}(z)\big) - \E \, \mathfrak{Re}\big(  S_{{{\bf B}_{n,m}}}(z)\big) \big| \!> 2x \Big) +\mathbb{P}\Big( \!\big| \mathfrak{Im}\big( S_{{{\bf B}_{n,m}}}(z)\big) - \E \, \mathfrak{Im}\big(  S_{{{\bf B}_{n,m}}}(z)\big) \big|\! > 2x \Big) 
\end{multline}
For a row vector $\textbf{x} \in \mathbb{R}^{N  n}$, we partition it into $n$ elements of dimension $N$ and write $\textbf{x} =\big (  x_{1}, \dots, x_{n} \big )$ where $ x_1, \ldots, x_{n}$ are row vectors of $\mathbb{R}^N $. Now, let $A(\textbf{x})$ and $B(\textbf{x})$ be respectively the $N \times n$ and $N \times N$ matrices defined by
\beq\label{deffA}
A(\textbf{x}) = \big(x_1^T | \ldots | x_n^T \big) 
\quad {\rm and} \quad
B(\textbf{x})=\frac{1}{n} A(\textbf{x})A(\textbf{x})^T.
\eeq
Also, let $h_1:= h_{1,z}$ and $h_2 := h_{2,z}$ be the functions defined from $\mathbb{R}^{Nn}$ into $\mathbb{R}$ by
$$
h_1({\bf x})= \int f_{1,z} \,  d \mu_{B({\bf x})}
\quad {\rm and} \quad
h_2({\bf x})= \int f_{2,z}  \, d \mu_{B({\bf x})} ,
$$
where $f_{1,z}(\lambda) =\frac{\lambda-u}{(\lambda-u)^2+v^2}$ and $f_{2,z}(\lambda) =\frac{v}{(\lambda-u)^2 + v^2}$ and note that
$
S_{B({\bf x})}(z) = h_1({\bf x})  + ih_2({\bf x}).
$
Now, denoting by ${\bf X}_{1,m}^T , \ldots , {\bf X}_{n,m}^T$ the columns of $\mathcal{X}_{n,m}$ and setting ${\bf A}$ to be the row random vector of $\mathbb{R}^{Nn}$ given by 
$$
{\bf A}= ( {\bf X}_{1,m} , \ldots , {\bf X}_{n,m}),
$$
we note that $B ({\bf A})={\bf B}_{n,m} $ and $ h_1({\bf A})= \mathfrak{Re}\big( S_{{{\bf B}_{n,m}}}(z)\big)$. Moreover, letting $q$ be a positive integer less than $n$, we set $\mathcal{F}_i= \sigma (\varepsilon_k, \ k \leqslant iN + m)$ for $1 \leqslant i \leqslant [n/q]q $ with the convention that $\mathcal{F}_0= \{ \emptyset , \Omega \}$ and that $\mathcal{F}_{s}=\mathcal{F}_n$ for any $s\in \{ [n/q]q, \ldots,n\} $. Noting that ${\bf X}_{1,m} , \ldots , {\bf X}_{i,m}$ are $\mathcal{F}_i$-measurable, we write the following decomposition:
\begin{align*}
\mathfrak{Re}\big( S_{{{\bf B}_{n,m}}}(z)\big) - \E \, \mathfrak{Re}\big(  S_{{{\bf B}_{n,m}}}(z)\big)
&=  h_1({\bf X}_{1,m} , \ldots , {\bf X}_{n,m}) - \E  h_1({\bf X}_{1,m} , \ldots , {\bf X}_{n,m}) \\
 &= \sum_{i=1}^{[n/q]} \big( \E ( h_1( {\bf A}) | \mathcal{F}_{iq})  - \E ( h_1( {\bf A}) | \mathcal{F}_{(i-1)q})\big). 
\end{align*}
Now, let $({\bf A}_i)_{ i }$ be a family of row random vectors of $ \mathbb{R}^{Nn}$ defined for any $i \in \{1 , \ldots , [n/q]-1\}$ by
$$
{\bf A}_i = \big({\bf X}_{1,m}\, , \ldots , {\bf X}_{(i-1)q,m}\, , \, \underbrace{{\bf 0}_{N}, \ldots  , {\bf 0}_{N}}_{2q \, {\rm times}}, {\bf X}_{(i+1)q+1,m} \, , \ldots , {\bf X}_{n,m} \big),
$$
and for $i=[n/q]$ by
$$
{\bf A}_{[\frac{n}{q}]} = \big({\bf X}_{1,m}\, , \ldots , {\bf X}_{([n/q]-1)q,m}\, , \, {\bf 0}_{N} , \ldots  , {\bf 0}_{N} \big).
$$
Noting that $\E \big(h_1({\bf A}_{[n/q]}) | \mathcal{F}_{n}\big)= \E (h_1 \big({\bf A}_{[n/q]}) | \mathcal{F}_{([n/q]-1)q}\big) $, we write
\begin{multline}\label{mart dcmp of Stiel 1}
\mathfrak{Re}\big( S_{{{\bf B}_{n,m}}}(z)\big)- \E \, \mathfrak{Re}\big(  S_{{{\bf B}_{n,m}}}(z)\big)\\
=\sum_{i=1}^{[n/q]} \Big( \E \big( h_1( {\bf A})- h_1( {\bf A}_i) | \mathcal{F}_{iq}\big)  - \E  \big( h_1( {\bf A})-h_1( {\bf A}_i)| \mathcal{F}_{(i-1)q}\big)\Big) \\
  +\sum_{i=1}^{[n/q]-1} \Big( \E \big( h_1( {\bf A}_i) | \mathcal{F}_{iq} \big) \!- \! \E \big( h_1( {\bf A}_i)| \mathcal{F}_{(i-1)q}\big)\Big)  \\
 := M_{[n/q],\, q} + \sum_{i=1}^{[n/q]-1} \Big( \E \big( h_1( {\bf A}_i) | \mathcal{F}_{iq}\big)  - \E \big( h_1( {\bf A}_i)| \mathcal{F}_{(i-1)q}\big)\Big).
\end{multline}
Thus, we get
\begin{align}\label{sum two Prob}
\mathbb{P}\Big( \Big| & \mathfrak{Re}\big( S_{{{\bf B}_{n,m}}}(z)\big) - \E \, \mathfrak{Re}\big(  S_{{{\bf B}_{n,m}}}(z)\big) \Big| > 2x \Big) \notag \\
&\leqslant \mathbb{P} \left( |M_{[n/q], \, q}  |>x \right)
  + \mathbb{P} \Big( \Big|\sum_{i=1}^{[n/q]-1 } \big( \E (h_1({\bf A}_i)| \mathcal{F}_{iq}) 
	- \E (h_1({\bf A}_i)| \mathcal{F}_{(i-1)q})\big) \Big|>x \Big).
\end{align}
Note that $( M_{k, q })_{k }$ is a centered martingale with respect to the filtration $(\mathcal{G}_{k,q})_{k}$ defined by $\mathcal{G}_{k,q}= \mathcal{F}_{kq}$. Moreover, for any $k \in \{ 1, \ldots , [n/q]\}$,
\begin{align*}
\Vert M_{k,q} - M_{k-1,q} \Vert_{\infty}
  &=  \Vert  \E ( h_1({\bf A})- h_1({\bf A}_k)| \mathcal{F}_{kq}) 
	  - \E ( h_1({\bf A})- h_1({\bf A}_k)| \mathcal{F}_{(k-1)q}) \Vert_{\infty} \\ 
	& \leqslant 2 \Vert h_1({\bf A})- h_1({\bf A}_k) \Vert_{\infty} 
\end{align*}
Noting that $\Vert f_{1,z}' \Vert_1=2/v$ then by integrating by parts, we get 
\begin{align}\label{Rank IPP Re(S)}
\vert h_1({\bf A})- h_1({\bf A}_k) \vert &= \Big| \int f_{1,z} d \mu_{B({\bf A})} -  \int f_{1,z} d \mu_{B({\bf A}_{k})} \Big| \leqslant \Vert f_{1,z}' \Vert_1 \Vert F^{B({\bf A})} - F^{B({\bf A}_k)} \Vert_{\infty}\notag\\
& \leqslant \frac{2}{v N} {\rm Rank} \big ( A({\bf A}) - A({\bf A}_k) \big),
\end{align}
where the second inequality follows from Theorem A.44 in \cite{BS}. As for any $k \in \{1 , \ldots , [n/q]-1\}$, $ {\rm Rank} \big ( A({\bf A}) - A({\bf A}_k) \big) \leqslant 2q$ and $ {\rm Rank} \big ( A({\bf A}) - A({\bf A}_{[n/q]}) \big) \leqslant q$, then overall we derive that almost surely
$$
\Vert M_{k,q} - M_{k-1,q} \Vert_{\infty} \leqslant \frac{8 q}{v N} \ \ \text{ and \ \ } 
\Vert M_{[n/q],q} - M_{[n/q]-1,q} \Vert_{\infty} \leqslant \frac{4 q}{v N} 
$$
and hence applying the Azuma-Hoeffding inequality for martingales we get for any $x>0$,
\begin{align}\label{conc ineq martingale part 1}
\mathbb{P} \left( |M_{[n/q],\, q}  |>x \right)
	\leqslant 2 \exp\left\lbrace - \frac{x^2 v^2  N^2}{128\, q \, n } \right\rbrace.
\end{align}
Now to control the second term of \eqref{mart dcmp of Stiel 1}, we have, by Markov's inequality and orthogonality, for any $x>0$,
\begin{multline}\label{proba second part conc ineq}
\mathbb{P} \Big( \Big|\sum_{i=1}^{[n/q]-1 } \big( \E (h_1({\bf A}_i)| \mathcal{F}_{iq}) 
	- \E (h_1({\bf A}_i)| \mathcal{F}_{(i-1)q})\big) \Big| >x \Big)
	\\
	 \leqslant \frac{1}{x^2} \sum_{i=1}^{[n/q]-1 } \Vert \E (h_1({\bf A}_i)| \mathcal{F}_{iq}) 
	- \E (h_1({\bf A}_i)| \mathcal{F}_{(i-1)q}) \Vert_2^2.
\end{multline} 
Fixing $i \in \{ 1, \ldots ,[n/q]-1\}$, one can construct by Berbee's maximal coupling lemma \cite{Berbee}, a sequence $(\varepsilon_k')_{k \in \mathbb{Z} }$ distributed as $(\varepsilon_k)_{k \in \mathbb{Z} }$ and independent of $ \mathcal{F}_{iq}$ such that for any $j > iqN+m$,
\begin{equation} \label{Proba from maximal coupling}
\mathbb{P} (  \varepsilon_k' \neq \varepsilon_k, \; {\rm for \; some} \; k \geqslant j ) = \beta_{j-iqN-m}.
\end{equation}
Let $({X}'_{k, m})_{k \geqslant 1 }$ be the sequence defined for any $k \geqslant 1$ by
$
{X}'_{k, m}= \E ( X_k | \varepsilon_{k-m}', \ldots , \varepsilon_{k+ m}')
$
and let $ {\bf X}_{i,m}'$ be the row vector of $\mathbb{R}^N$ defined by
${\bf X}_{i,m}'=({ X}_{(i-1)N+1, m}', \ldots , { X}_{iN, m}')$. Finally, we define for any $i \in \{1 , \ldots , [n/q]-1\}$ the row random vector ${\bf A}_i'$ of $ \mathbb{R}^{Nn}$ by
$$
{\bf A}_i' = \big({\bf X}_{1,m}\, , \ldots , {\bf X}_{(i-1)q,m}\, , \underbrace{ {\bf 0}_{N}, \ldots  , {\bf 0}_{N}}_{2q \, {\rm times}}, \,  {\bf X}_{(i+1)q+1,m}' \, , \ldots , {\bf X}_{n,m}'\, \big).
$$
As $ {\bf X}_{(i+1)q+1,m}'  , \ldots , {\bf X}_{n,m}' $ are independent of $\mathcal{F}_{iq}$ then $\E(h_1({\bf A}_i')|\mathcal{F}_{iq})= \E(h_1({\bf A}_i')|\mathcal{F}_{(i-1)q})$. Thus we write
$$
\E\big( h_1({\bf A}_i) |\mathcal{F}_{iq}\big) - \E\big( h_1({\bf A}_i) |\mathcal{F}_{(i-1)q}\big) 
= \E\big( h_1({\bf A}_i) - h_1({\bf A}_i') |\mathcal{F}_{iq}\big) - \E\big( h_1({\bf A}_i) - h_1({\bf A}_i')|\mathcal{F}_{(i-1)q}\big).
$$
and infer that
\begin{multline} \label{conc ineq UB 2}
\Vert \E\big( h_1({\bf A}_i) |\mathcal{F}_{iq}\big)  - \E\big( h_1({\bf A}_i) |\mathcal{F}_{(i-1)q}\big)\Vert_2  \\
	 \leqslant \Vert \E\big( h_1({\bf A}_i) - h_1({\bf A}_i') |\mathcal{F}_{iq}\big) \Vert_2
	 + \Vert \E\big( h_1({\bf A}_i) - h_1({\bf A}_i')|\mathcal{F}_{(i-1)q}\big)\Vert_2 \\
	 \leqslant  2 \Vert h_1({\bf A}_i) - h_1( {\bf A}_{i}')  \Vert_2
\end{multline} 
Similarly as in \eqref{Rank IPP Re(S)}, we have
\begin{align*}
\big| h_1({\bf A}_i) - h_1( {\bf A}_i')  \big| 
& \leqslant \frac{2}{v N} {\rm Rank} (A({\bf A}_i) -  A({\bf A}_i')) 
 \leqslant \frac{2}{v N} \sum_{\ell= (i+1)q+1}^{n}  {\bf 1}_{ \big\lbrace {\bf X}_{\ell,m }' \neq {\bf X}_{\ell,m}  \big\rbrace}  \\
& \leqslant \frac{2 n }{v N}\,  {\bf 1}_{ \big\lbrace \varepsilon_k' \neq \varepsilon_k, \ {\rm for \ some \; } k \geqslant \, (i+1)qN +1-m \big\rbrace }\, .
\end{align*}
Hence by \eqref{Proba from maximal coupling}, we infer that 
\begin{align}\label{conc beta}
\Vert h_1({\bf A}_i) - h_1( {\bf A}_{i}')  \Vert_2^2 &\leqslant \frac{4 n^2 }{v^2 N^2} \,\beta_{qN +1-2m} 
\leqslant \frac{4 n^2 }{v^2 N^2}\,  \beta_{(q-1)N} ,
\end{align}
Starting from \eqref{proba second part conc ineq} together with  \eqref{conc ineq UB 2} and \eqref{conc beta}, it follows that 
\beq\label{proba coupling meth}
\mathbb{P} \Big( \Big|\sum_{i=1}^{[n/q]-1 } \big( \E (h_1({\bf A}_i)| \mathcal{F}_{iq}) 
	- \E (h_1({\bf A}_i)| \mathcal{F}_{(i-1)q})\big) \Big|>x \Big) 
	\leqslant\frac{16 n^3\, }{x^2 \, v^2 q \, N^2} \beta_{(q-1)N}.
\eeq
Therefore, considering \eqref{sum two Prob} and gathering the upper bounds in \eqref{conc ineq martingale part 1} and \eqref{proba coupling meth}, we get 
$$
\mathbb{P} \big( \big| \mathfrak{Re}\big( S_{{{\bf B}_{n,m}}}(z)\big) - \E \, \mathfrak{Re}\big(  S_{{{\bf B}_{n,m}}}(z)\big)\big| > 2x \big) \leqslant 2 \exp\left\lbrace - \frac{x^2 v^2  N^2}{128\, q \, n }\right\rbrace 
+\frac{16 n^3\, }{x^2 \, v^2 q \, N^2} \beta_{(q-1)N}.
$$
Finally, noting that $ \mathbb{P} \big( \big| \mathfrak{Im}\big( S_{{{\bf B}_{n,m}}}(z)\big) - \E \, \mathfrak{Im}\big(  S_{{{\bf B}_{n,m}}}(z)\big)\big| > 2x \big)$ also admits the same upper bound and choosing $q= [n/(\log n)^{\alpha} ] +1 $, \eqref{conc ineq stieltjes Bn,m bar} follows. This ends the proof of Theorem \ref{theo conc ineq}.

\section{Proof of Theorem \ref{first theorem}}
\setcounter{equation}{0}
The proof, being technical, will be divided into three major steps (Sections \ref{app Bn by Bn,m} to \ref{section construction of gaussians}).

\subsection{A first approximation} \label{app Bn by Bn,m}
 Let $m$ be a fixed positive integer and set $p:= p(m)=a_m m$ with $(a_n)_{n \geqslant 1}$ being a sequence of positive integers such that $\lim_{n\rightarrow \infty} a_n=\infty$.
Setting $k_N= \Big[ \frac{N}{p+3m} \Big]$, we write the subset $ \lbrace 1 , \ldots , Nn \rbrace$ as a union of disjoint subsets of $\mathbb{N}$ as follows:
$$
[1, Nn] \cap \mathbb{N}= \bigcup_{i=1}^{n} [(i-1)N+1, iN] \cap \mathbb{N} =
 \bigcup_{i=1}^{n} \bigcup_{\ell=1}^{k_N +1} I^{i}_{\ell} \cup J^{i}_{\ell}, 
$$
where, for $i  \in \{ 1, \ldots , n\} $ and $ \ell \in \{1, \dots, k_N \}$,
\begin{align*}
 I^{i}_{\ell} &:=\big [(i-1)N + (\ell-1)(p+3m)+1  \, , \, (i-1)N + (\ell-1)(p+ 3m )+p \big ]\cap \mathbb{N} , \label{defI}\\
 J^{i}_{\ell} &:= \Big [ (i-1)N + (\ell-1)(p+3m)+p+1 \, , \, (i-1)N + \ell(p+3m)\Big ]\cap \mathbb{N}  \, ,
\end{align*}
and, for $\ell = k_{N} +1 $, $I^{i}_{k_{N} +1}  = \emptyset $ and
$$
J^{i}_{k_{N} +1}  = \big [(i-1)N + k_{N}(p+ 3m)+1 \, , \, iN \big ] \cap \mathbb{N} \, .
$$
Note that for all $ i \in \{ 1, \ldots , n \}$, $J^{i}_{k_{N} +1}  =\emptyset$ if $k_{N}(p+3m) = N$.
Now, let $M$ be a fixed positive number not depending on $(n,m)$ and  let $\varphi_M$ be the function defined by $\varphi_M(x)= (x \wedge M) \vee (-M)$.  Setting
\beq\label{Bi,l def}
B_{i, \ell} = ( \varepsilon_{(i-1)N + (\ell-1)(p+3m)+1-m}, \ldots , \varepsilon_{(i-1)N + (\ell-1)(p+3m)+p+ m} ),
\eeq
we define the sequences $(\widetilde{X}_{k,m,M})_{k \geqslant 1 }$ and $(\bar{X}_{k,m,M})_{k \geqslant 1 }$ as follows:
\[
\widetilde{X}_{k, m,M}= \begin{cases}
\E (\varphi_M( X_k) | B_{i, \ell}) & \textrm{if  \; $k \in   I^i_{\ell} \;$} \\
0  & \textrm{otherwise,}\\
\end{cases}
\]
and 
\begin{equation}\label{def of Xk,m}
\bar{X}_{k,m,M}= \widetilde{X}_{k, m,M} - \E (\widetilde{X}_{k, m,M}).
\end{equation}
To soothe the notations, we shall write $\widetilde{X}_{k, m}$ and $\bar{X}_{k,m} $ instead of $\widetilde{X}_{k, m,M}$ and $\bar{X}_{k,m,M}$ respectively. Note that for any $k \geqslant 1 $,
\begin{equation} \label{norm2 Xk,m}
\Vert \bar{X}_{k,m} \Vert_2 \leqslant 2 \Vert \widetilde{X}_{k,m} \Vert_2 = 2 \Vert \E (\varphi_M(X_k)|B_{i,\ell}) \Vert_2 \leqslant 2 \Vert \varphi_M(X_k) \Vert_2  \leqslant 2 \Vert X_k \Vert_2 = 2\Vert X_0 \Vert_2,
\end{equation}
and
\begin{equation}\label{norm infty Xk,m}
\Vert \bar{X}_{k,m} \Vert_{\infty} \leqslant 2  \Vert  \widetilde{X}_{k,m} \Vert_{\infty} \leqslant 2 M ,
\end{equation}
where the last equality in \eqref{norm2 Xk,m} follows from the stationarity of $(X_k)_{k }$.  As $\bar{X}_{k,m}$ is $\sigma (B_{i, \ell})$-measurable then it can be written as a measurable function $h_k$ of $B_{i, \ell}$, i.e. 
\begin{equation}\label{X-k,m-function-hk}
\bar{X}_{k,m} =h_k ( B_{i, \ell} ) .
\end{equation}
Finally, let $\bar{\mathcal{X}}_{n,m}= ((\bar{\mathcal{X}}_{n,m})_{i,j})=(\bar{X}_{(j-1)N +i ,\, m})$ and set
\begin{equation} \label{matrix Bn,m}
\bar{\bf B}_{n,m} = \frac{1}{n} \bar{\mathcal{X}}_{n,m} \bar{\mathcal{X}}_{n,m}^T \, .
\end{equation}
We shall approximate ${\bf B}_n$ by $\bar{\bf B}_{n,m}$ by applying the following proposition:
\begin{prop}\label{prop Bn Bn,m}
Let ${\bf B}_n$ and $\bar{\bf B}_{n,m}$ be the matrices defined in \eqref{matrix-Bn} and \eqref{matrix Bn,m} respectively then if $\lim_{n\rightarrow \infty} N/n = c \in (0, \infty)$, we have for any $z \in \mathbb{C}^{+}$,
\begin{equation}\label{app-Btilde-by-B_n,m}
\lim_{m\rightarrow + \infty} \limsup_{M \rightarrow + \infty}  \limsup_{n \rightarrow + \infty} \big| \E \big(S_{{{\bf B}_n}}(z)\big) -\E \big(S_{{\bar{\bf B}_{n,m}}}(z)\big)\big| = 0 .
\end{equation}
\end{prop}

\noindent
\textit{Proof.} By \eqref{BaMe prop 4.1} and Cauchy-Schwarz's inequality, it follows that
\begin{align}\label{trace lemma Bn Bn,m}
\big| \E \big(S_{{{\bf B}_n}}(z)\big) &-\E \big(S_{{\bar{\bf B}_{n,m}}}(z)\big)\big| \notag \\
&\leqslant \frac{\sqrt{2}}{v^2} \left\Vert \frac{1}{N} {\rm Tr}({\bf B}_n + \bar{\bf B}_{n,m} ) \right\Vert_1^{1/2}
		  \left\Vert \frac{1}{Nn} {\rm Tr} ( \mathcal{X}_n - \bar{\mathcal{X}}_{n,m}) ( \mathcal{X}_n - \bar{\mathcal{X}}_{n,m})^T \right\Vert_1^{1/2}\, .
\end{align}
By the definition of ${\bf B}_n$, $ N^{-1}\E \big| {\rm Tr}({\bf B}_n)\big|= \Vert X_0 \Vert_2^2$. Similarly and due to the fact that $p k_N \leqslant N$ and \eqref{norm2 Xk,m}, 
\beq \label{trace Bn,m}
 \frac{1}{N}\E \big| {\rm Tr}(\bar{\bf B}_{n,m})\big|= \frac{1}{Nn} \sum_{i=1}^n \sum_{\ell=1}^{k_N} \sum_{k \in I^i_{\ell}}  \Vert \bar{X}_{k,m}\Vert_2^2 \leqslant 4 \Vert X_0 \Vert_2^2 \, .
\eeq
Moreover, by the construction of $\mathcal{X}_n$ and $\bar{\mathcal{X}}_{n,m}$, we have
$$
\frac{1}{Nn}  \E \vert {\rm Tr} ( \mathcal{X}_n - \bar{\mathcal{X}}_{n,m}) ( \mathcal{X}_n - \bar{\mathcal{X}}_{n,m})^T \vert = \frac{1}{Nn} \sum_{i=1}^n \sum_{\ell=1}^{k_N} \sum_{k \in I^i_{\ell}}  \Vert X_k - \bar{X}_{k,m}\Vert_2^2 + \frac{1}{Nn} \sum_{i=1}^n \sum_{\ell=1}^{k_N+1} \sum_{k \in J^i_{\ell}}  \Vert  X_{k}\Vert_2^2 \, .
$$
Now, since $X_k$ is centered, we write for $k \in I_{\ell}^i$,  
\begin{align}\label{Xk-Xk,mbar-norme2}
\Vert X_k - \bar{X}_{k,m}\Vert_2 
	&= \Vert X_k  - \widetilde{X}_{k,m} -\E(X_k  - \widetilde{X}_{k,m})  \Vert_2 
	\leqslant  2 \Vert X_k  - \widetilde{X}_{k,m} \Vert_2\notag \\
	& \leqslant 2 \Vert X_k  - \E (X_k|B_{i,\ell}) \Vert_2 +  2\Vert \widetilde{X}_{k,m} -\E (X_k|B_{i,\ell}) \Vert_2\, .
\end{align}
Analyzing the second term of the last inequality, we get
\begin{equation}\label{X - M+}
\Vert \widetilde{X}_{k,m} -\E (X_k|B_{i,\ell}) \Vert_2 
= \Vert \E( X_{k} - \varphi_M(X_k)|B_{i,\ell}) \Vert_2
\leqslant \Vert X_{k} - \varphi_M(X_k )\Vert_2
= \Vert ( |X_0| - M)_{+} \Vert_2.
\end{equation}
As $X_0$ belongs to $\mathbb{L}^2$, then $\lim_{M \rightarrow + \infty} \Vert ( |X_{0}| - M)_{+} \Vert_2 =0$. Now, we note that for $ k \in I^i_{\ell} $, $ \sigma ( \varepsilon_{k-m} , \ldots , \varepsilon_{k+m}) \subset  \sigma (B_{i,\ell})$ which implies that
\begin{align}\label{xk - Xk,m hat}
\Vert X_k  - \E (X_k|B_{i,\ell}) \Vert_2
&\leqslant \Vert X_k - \E( X_k| \varepsilon_{k-m} , \ldots , \varepsilon_{k+m}) \Vert_2 \notag\\
& = \Vert X_0 -\E( X_0| \varepsilon_{-m} , \ldots , \varepsilon_{m}) \Vert_2 
=\Vert X_0 -{X}_{0,m}\Vert_2\, ,
\end{align}
where the first equality is due to the stationarity. Therefore, by \eqref{X - M+}, \eqref{xk - Xk,m hat}, the fact that $pk_N \leqslant N$ and  
\[
{\rm Card} \Big( \bigcup_{i=1}^n \bigcup_{\ell = 1}^{k_N+1} J_{\ell}^{i} \Big) \leqslant Nn - npk_N\, ,
\] 
we infer that
\begin{align}\label{bound diff tr Xn Xn,m}
\frac{1}{Nn}  \E \vert {\rm Tr} ( \mathcal{X}_n - \bar{\mathcal{X}}_{n,m}) ( \mathcal{X}_n - \bar{\mathcal{X}}_{n,m})^T \vert
&\leqslant 8\Vert X_0 -{X}_{0,m}\Vert_2^2 
+ 8\Vert ( |X_{0}| - M)_{+} \Vert_2^2 \notag \\
& \quad \quad+(3(a_m +3)^{-1}+a_m m N^{-1})\Vert X_0 \Vert_2^2.
\end{align}
Thus starting from \eqref{trace lemma Bn Bn,m}, considering the upper bounds \eqref{trace Bn,m} and \eqref{bound diff tr Xn Xn,m}, we derive that there exists a positive constant $C$ not depending on $(n,m,M)$ such that
$$
\limsup_{M\rightarrow \infty}\limsup_{n \rightarrow \infty} \big| \E \big(S_{{{\bf B}_n}}(z)\big)-\E \big(S_{{\bar{\bf B}_{n,m}}}(z)\big)\big|
\leqslant\frac{C}{v^2} \Big( \Vert X_0 -{X}_{0,m}\Vert_2^2 +\frac{3}{a_m}\Big)^{1/2}.
$$
Taking the limit on $m$, Proposition \ref{prop Bn Bn,m} follows by applying the martingale convergence theorem in $\mathbb{L}^2$ and that fact that $a_m$ converges to infinity. 
\hfill$\square$

\subsection{Approximation by a Gram matrix with independent blocks}

By Berbee's classical coupling lemma \cite{Berbee}, one can construct by induction a sequence of random variables $(\varepsilon_{k}^{*})_{k \geqslant 1}$ such that:
\begin{itemize}
\item For any $1 \leqslant i\leqslant n$ and $1 \leqslant \ell \leqslant k_N$, $$
B_{i , \ell}^{*}=(\varepsilon_{(i-1)N+(\ell -1)(p+3m)+1-m}^{*}, \ldots , \varepsilon_{(i-1)N+(\ell -1)(p+3m)+p+m}^{*} )
$$
has the same distribution as $B_{i , \ell}$ defined in \eqref{Bi,l def}.
\item The array $(B_{i , \ell}^{*})_{1 \leqslant i\leqslant n ,\, 1 \leqslant \ell \leqslant k_N}$ is i.i.d.
\item For any $1 \leqslant i\leqslant n$ and $1 \leqslant \ell \leqslant k_N$, $\mathbb{P}(B_{i, \ell} \neq B_{i, \ell}^{*}) \leqslant \beta_m$. 
\end{itemize}
(see page 484 of \cite{Viennet} for more details concerning the construction of the array $(B_{i , \ell}^{*})_{i, \ell \geqslant 1}$). We define now the sequence $ ( \bar{X}_{k,m}^{*})_{k \geqslant 1}$ as follows:
\begin{equation} \label{defi of Xk,m star}
\bar{X}_{k,m}^{*}= h_k( B_{i, \ell}^{*}) \quad {\rm if } \; k \in   I^i_{\ell}, 
\end{equation}
where the functions $h_k$ are defined in (\ref{X-k,m-function-hk}).

\smallskip
\noindent
We construct the $N \times n$ random matrix $\bar{\mathcal{X}}_{n,m}^{*}= ((\bar{\mathcal{X}}_{n,m}^{*})_{i,j})=(\bar{X}_{(j-1)N+i, \, m}^{*})$. Note that the block of entries $(\bar{X}_{k,m}^{*},\, k \in I^i_{\ell})$ is independent of $(\bar{X}_{k,m}^{*},\, k \in I^{i'}_{\ell '})$ if $(i,\ell)\neq (i' , \ell ')$. Thus, $\bar{\mathcal{X}}_{n,m}^{*}$ has independent blocks of dimension $p$ separated by null blocks whose dimension is at least $3m$.
Setting
\begin{equation}\label{matrix Bn,m Star}
\bar{\bf B}_{n,m}^{*} := \frac{1}{n}
\bar{\mathcal{X}}_{n,m}^{*} \bar{\mathcal{X}}_{n,m}^{* \,T}\, ,
\end{equation}
we approximate $\bar{\bf B}_{n,m}$ by the Gram matrix $\bar{\bf B}_{n,m}^{*}$ as shown in the following proposition. 
\begin{prop}\label{prop Bn,m Bn,m star}
Let $\bar{\bf B}_{n,m}$ and $\bar{\bf B}_{n,m}^{*}$ be defined in \eqref{matrix Bn,m} and \eqref{matrix Bn,m Star} respectively. Assuming that $\ \lim_{n\rightarrow \infty} N/n = c \in (0, \infty)$ and $\lim_{n \rightarrow \infty}\beta_n=0$ then for any $z \in \mathbb{C}^{+}$,
\begin{equation}\label{app-Bn,m-by-Bn,m tilde}
\lim_{m \rightarrow + \infty} \limsup_{M \rightarrow + \infty} \limsup_{n \rightarrow + \infty} \big| \E \big(S_{{\bar{\bf B}_{n,m}}}(z) \big) - \E \big(S_{{\bar{\bf B}_{n,m}^{*}}}(z) \big) \big| = 0 .
\end{equation}
\end{prop}

\noindent
\textit{Proof.} By \eqref{BaMe prop 4.1} and Cauchy-Schwarz's inequality, it follows that
\begin{multline}\label{trace lemma Bn,m Bn,m star}
\big| \E \big(S_{{\bar{\bf B}_{n,m}}}(z)\big) -\E \big(S_{{\bar{\bf B}_{n,m}^{*}}}(z)\big)\big|\\
\leqslant \frac{\sqrt{2}}{v^2} \left\Vert \frac{1}{N} {\rm Tr}(\bar{\bf B}_{n,m} + \bar{\bf B}_{n,m}^{*} ) \right\Vert_1^{1/2}
		  \left\Vert \frac{1}{Nn} {\rm Tr} ( \bar{\mathcal{X}}_{n,m} - \bar{\mathcal{X}}_{n,m}^{*}) ( \bar{\mathcal{X}}_{n,m} - \bar{\mathcal{X}}_{n,m}^{*})^T \right\Vert_1^{1/2}
\end{multline}
Notice that $ \Vert \bar{X}_{k,m}^{*} \Vert_2 = \Vert h_k ( B_{i , \ell}^{*}) \Vert_2 = \Vert h_k ( B_{i , \ell}) \Vert_2 
= \Vert \bar{X}_{k,m} \Vert_2 \leqslant  2 \Vert X_0 \Vert_2$, where the second equality follows from the fact that $ B_{i , \ell}^{*}$ is distributed as $ B_{i , \ell}$ whereas the last inequality follows from (\ref{norm2 Xk,m}). Thus, we get from the definition of $\bar{\bf B}_{n,m}^{*}$ and the fact that $pk_N \leqslant N$,
\begin{align}\label{tr Bn,m star}
 \frac{1}{N} \E \big| {\rm Tr} (  \bar{\bf B}_{n,m}^{*})\big|
 	  =  \frac{1}{N n } \sum_{i=1}^{n} \sum_{\ell = 1}^{k_N} \sum_{k \in I^{i}_{\ell} } \Vert \bar{X}_{k,m}^{*}\Vert_2^2
   	   	\leqslant 4  \Vert X_0 \Vert_2^2.
\end{align}
Considering \eqref{trace lemma Bn,m Bn,m star}, \eqref{trace Bn,m} and \eqref{tr Bn,m star}, we infer that Proposition \ref{prop Bn,m Bn,m star} follows once we prove that
\beq\label{le but ds prop star}
\lim_{m \rightarrow + \infty} \limsup_{M \rightarrow + \infty} \limsup_{n \rightarrow + \infty} \frac{1}{Nn} \E \Big|{\rm Tr} (\bar{\mathcal{X}}_{n,m} -\bar{\mathcal{X}}_{n,m}^{*})( \bar{\mathcal{X}}_{n,m} - \bar{\mathcal{X}}_{n,m}^{*})^T \Big|=0 \, .
\eeq

\noindent
By the construction of $\bar{\mathcal{X}}_{n,m}$ and $\bar{\mathcal{X}}_{n,m}^{*}$, we write
\begin{align}\label{Tr_Xn,m minus Xn,m tilde}
\frac{1}{Nn} \E \Big|{\rm Tr} (\bar{\mathcal{X}}_{n,m} &-\bar{\mathcal{X}}_{n,m}^{*})( \bar{\mathcal{X}}_{n,m} - \bar{\mathcal{X}}_{n,m}^{*})^T \Big|
	 = \frac{1}{Nn} \sum_{i =1}^{n} \sum_{\ell = 1}^{k_N} \sum_{k \in I^{i}_{\ell} } 
		\Vert \bar{X}_{k,m} - \bar{X}_{k,m}^{*}\Vert_2^2 .
\end{align}
Now, let $L$ be a fixed positive real number strictly less than $M$ and not depending on $(n,m,M)$. To control the term $\Vert \bar{X}_{k,m} - \bar{X}_{k,m}^{*} \Vert_2^2$, we write for  $k \in I_{\ell}^{i}$,
\begin{multline}
\Vert \bar{X}_{k,m} - \bar{X}_{k,m}^{*} \Vert_2^2
   = \Vert (h_k ( B_{i , \ell}) - h_k ( B_{i , \ell}^{*})) {\bf 1}_{B_{i , \ell} \neq B_{i , \ell}^{*}}\Vert_2^2\\
   \leqslant 4 \Vert h_k ( B_{i , \ell}) {\bf 1}_{B_{i , \ell} \neq B_{i , \ell}^{*}}\Vert_2^2 
   = 4 \Vert \bar{X}_{k,m} {\bf 1}_{B_{i , \ell} \neq B_{i , \ell}^{*}}\Vert_2^2 \quad  \quad \quad \quad \quad \quad \quad \quad \quad\\
   \leqslant 12 \Vert \bar{X}_{k,m} - \E( X_k| B_{i , \ell})\Vert_2^2
    + 12  \Vert \E( X_k| B_{i , \ell}) - \E (\varphi_L (X_k)| B_{i , \ell}) \Vert_2^2 \\
     \quad \quad +12 \Vert \E (\varphi_L (X_k)| B_{i , \ell}) {\bf 1}_{B_{i , \ell} \neq B_{i , \ell}^{*}} \Vert_2^2 \, . \notag
 \end{multline}
Since $\mathbb{P}(B_{i, \ell} \neq B_{i, \ell}^{*}) \leqslant \beta_m$ and $\varphi_L(X_k)$ is bounded by $L$, we get
$$
\Vert \E (\varphi_L (X_k)| B_{i , \ell}) {\bf 1}_{B_{i , \ell} \neq B_{i , \ell}^{*}} \Vert_2^2 \leqslant L^2 \beta_m \,.
$$
Moreover, it follows from the fact that $X_k$ is centered and \eqref{X - M+} that
\[
 \Vert \bar{X}_{k,m} - \E( X_k| B_{i , \ell}) \Vert_2^2 
 \leqslant 4  \Vert \widetilde{X}_{k,m} - \E( X_k| B_{i , \ell}) \Vert_2^2
 \leqslant 4 \Vert (|X_0|-M)_{+} \Vert_2^2
\]
and 
\[
\Vert \E( X_k| B_{i , \ell}) - \E (\varphi_L (X_k)| B_{i , \ell}) \Vert_2^2
\leqslant \Vert  X_k - \varphi_L (X_k) \Vert_2^2
=\Vert (|X_0|-L)_{+} \Vert_2^2.
\]
Hence gathering the above upper bounds we get
\begin{equation}\label{norm Xbar - Xbar*} 
		\Vert \bar{X}_{k,m} - \bar{X}_{k,m}^{*}\Vert_2^2 
		\leqslant 48 \Vert (|X_0|-M)_{+} \Vert_2^2 + 12  \Vert (|X_0|-L)_{+} \Vert_2^2 + 12 L^2 \beta_m.
\end{equation}
As $p k_N \leqslant N$, the right-hand-side of \eqref{Tr_Xn,m minus Xn,m tilde} converges to zero by letting first $M $, then $m$ and finally $L$ tend to infinity.  Therefore, \eqref{le but ds prop star} and thus the proposition  follow.
\hfill$\square$

\subsection{Approximation with a Gaussian matrix} \label{section construction of gaussians}

\smallskip
\noindent

In order to complete the proof of the theorem, it suffices, in view of \eqref{app-Btilde-by-B_n,m} and \eqref{app-Bn,m-by-Bn,m tilde}, to prove the following convergence: for any $z \in \mathbb{C}^{+}$,
\begin{equation}\label{app by gaussian process}
\lim_{m\rightarrow + \infty} \limsup_{n \rightarrow + \infty} \big| \E \big ( S_{{{\bar{\bf B}_{n,m}^{*}}}}(z) \big )  -\E \big ( S_{{{\bf{G}}_{n}}}(z) \big ) \big| =0 .
\end{equation}
With this aim, we shall first consider a sequence $(Z_{k,m})_{k\in \mathbb{Z}}$ such that for any $k, \ell \in \{ 1, \ldots, N\}$,
\begin{equation}\label{cov-Zk,m}
\text{Cov}(Z_{k,m} \,, Z_{\ell,m})
=\text{Cov} ( \bar{X}_{k,m}^{*} \,,\bar{X}_{\ell,m}^{*})
\end{equation}
and let $(Z_{k,m}^{(i)})_k$, $i =1, \ldots ,n$, be $n$ independent copies of $(Z_{k,m})_k$. 
We then define the $N \times n$ matrix $\mathcal{Z}_{n,m} = ((\mathcal{Z}_{n,m})_{u,v})= (Z_{u,m}^{(v)})$ and finally set
\beq 
{\bf{G}}_{n,m}= \frac{1}{n} \mathcal{Z}_{n,m} \mathcal{Z}_{n,m}^{T} \, .
\eeq 

Now, we shall construct a matrix $\widetilde{\mathcal{Z}}_{n,m}$ having the same block structure as the matrix $\bar{\mathcal{X}}_{n,m}^{*}$. With this aim, we let for $\ell=1, \ldots ,k_N$,
$$
I_{\ell}=\{ (\ell-1)(p+3m)+1, \ldots , (\ell-1)(p+3m)+p \}
$$
and let $\widetilde{Z}_{k,m}^{(i)}$ be defined for any $ 1 \leqslant i \leqslant n$ and $ 1 \leqslant k \leqslant N$ by
\[
\widetilde{Z}_{k,m}^{(i)}= \begin{cases}
Z_{k,m}^{(i)} & \textrm{if  \; $ k \in   I_{\ell} \;$ for  some  $ \ell \in \{  1, \ldots ,  k_N \} $} \\
0  & \textrm{otherwise\, .}\\
\end{cases}
\]
We define now the $N \times n$ matrix $\widetilde{\mathcal{Z}}_{n,m} = ((\widetilde{\mathcal{Z}}_{n,m})_{u,v})= (\widetilde{Z}_{u,m}^{(v)})$ and we note that $\widetilde{\mathcal{Z}}_{n,m}$, as $\bar{\mathcal{X}}_{n,m}^{*}$, consists of independent blocks of dimension $p$ separated by null blocks whose dimension is at least $3m$. We finally set
\beq \label{def Gn bar}
\widetilde{\bf{G}}_{n,m}= \frac{1}{n} \widetilde{\mathcal{Z}}_{n,m}  \widetilde{\mathcal{Z}}_{n,m}^{T} \, .
\eeq 
Provided that $ \lim_{n\rightarrow \infty} n/N = c \in (0, \infty)$, we have by Proposition 4.2 in \cite{BaMe} that for any $z\in \mathbb{C}^+$,
\[ 
\lim_{m \rightarrow \infty}  \limsup_{n \rightarrow \infty}\Big | \E \big ( S_{{\bf G}_{n,m}}(z)  \big ) -\E \big ( S_{{\widetilde{\bf{G}}_{n,m}}}(z) \big ) \Big | = 0.
\]
In order to prove \eqref{app by gaussian process}, we shall prove for any $z\in \mathbb{C}^+$,
\[
\lim_{n \rightarrow + \infty} \big| \E \big ( S_{{{\bar{\bf B}_{n,m}^{*}}}}(z) \big )  -\E \big ( S_{{\widetilde{\bf{G}}_{n,m}}}(z) \big ) \big| =0 
\]
and then 
\[
\lim_{m \rightarrow \infty} \limsup_{M\rightarrow\infty} \limsup_{n \rightarrow \infty}\Big | \E \big ( S_{{{\bf G}_n}}(z)  \big ) -\E \big ( S_{{\bf{G}}_{n,m}}(z) \big ) \Big | = 0 \, .
\]
The technique followed to prove the first convergence is based on Lindeberg's method by blocks, whereas, for the second, it is based on the Gaussian interpolation technique.

\begin{prop}
Provided that $N/n\rightarrow c \in (0, \infty)$, then for any $z=x+iy \in \mathbb{C}^{+}$,
\begin{equation}\label{last step}
 \lim_{n \rightarrow + \infty} \big| \E \big ( S_{{{\bar{\bf B}_{n,m}^{*}}}}(z) \big )  -\E \big ( S_{{\widetilde{\bf{G}}_{n,m}}}(z) \big ) \big| =0 .
\end{equation}
\end{prop}

We don't give a full proof of this proposition because the computation involved has been almost exactly done in the proof of Proposition 4.3 of \cite{BaMe}. Indeed as $\bar{\mathcal{X}}_{n,m}^{*}$, the matrix $\bar{\mathcal{X}}_n$ considered in \cite{BaMe} has independent blocks separated by blocks of zero entries. The main difference is that in our case even the Gaussian blocks are mutually independent and thus the terms of first and second order in the Taylor expansion vanish simplifying the proof even more. Therefore a careful analysis of the proof of Proposition 4.3 in \cite{BaMe} gives for any $z=x+iy \in \mathbb{C}_{+}$,
$$
|\E\big (S_{\bar{\mathbf{B}}_{n,m}^{*}}(z)\big )
-\E \big (S_{\widetilde{\mathbf{G}}_{n,m}}(z)\big )| \leqslant \frac{Cp^2(1+M^3)N^{1/2}}{y^3 ( 1 \wedge y)n}  
$$
which converges to $0$ as $n$ tends to infinity since $\lim_{n \rightarrow \infty} N/n= c \in (0, \infty)$.
\hfill$\square$

\smallskip

To end the proof of Theorem \ref{first theorem}, it remains to prove the following convergence.

\begin{prop}\label{prop-gau-interpolation}
Provided that $ \lim_{n\rightarrow \infty} N/n = c \in (0, \infty)$ and $ \lim_{n\rightarrow \infty} \beta_n=0$ then for any $z\in \mathbb{C}^+$,
\[ 
\lim_{m \rightarrow \infty} \limsup_{M\rightarrow\infty} \limsup_{n \rightarrow \infty}\Big | \E \big ( S_{{{\bf G}_n}}(z)  \big ) -\E \big ( S_{{\bf{G}}_{n,m}}(z) \big ) \Big | = 0.
\]
\end{prop}

\noindent
\textit{Proof.} Let $n'=N+n$ and $\mathbb{G}_{n'}$ the symmetric matrix of order $n$ defined by
\[
\mathbb{G}_{n'}=\frac{1}{\sqrt{n}}\left(
\begin{array}
[c]{cc}%
\mathbf{0} & {\mathcal{Z}}_{n}^{T}\\
{\mathcal{Z}}_{n} & \mathbf{0} %
\end{array}
\right)  \,.
\]
It is well known that for any $z\in{\mathbb{C}}^{+}$,
\[
S_{\mathbf{G}_{n}}(z)=z^{-1/2}\frac{n}{2N}S_{\mathbb{G}_{n'}}(z^{1/2}%
)+\frac{n-N}{2Nz}\, 
\]
(See, for instance, page 549 in Rashidi Far \textit{et al} \cite{ROBS} for
arguments leading to the relation above). Since the same relation also holds
for the symmetric matrix $\mathbb{G}_{n',m}$ associated with $\mathbf{G}_{n,m}$ and
since $n'/N\rightarrow1+c^{-1}$, it is equivalent 
to prove for any $z\in{\mathbb{C}}^{+}$,
\[
\lim_{m\rightarrow\infty}\limsup_{M\rightarrow\infty}\limsup_{n\rightarrow\infty}\big |\E \big (S_{{\mathbb{G}}_{n'}}(z)\big)-\E \big (S_{\mathbb{G}_{n',m}}(z)\big )\big |=0\, .
\]
This can be proved in the same way as the convergence (32) in \cite{BaMePe}. Indeed, noting that $(\mathbb{G}_{n'})_{k,\ell}= n^{-1/2} Z_{k-n}^{(\ell)} {\bf 1}_{k >n} {\bf 1}_{ \ell \leqslant n}$  and $(\mathbb{G}_{n',m})_{k,\ell}= n^{-1/2} Z_{k-n,m}^{(\ell)} {\bf 1}_{k >n} {\bf 1}_{ \ell \leqslant n}$ if $1 \leqslant \ell \leqslant k \leqslant n' $ and keeping in mind the independence structure between the columns, we apply Lemma 16 in \cite{BaMePe} and we get for any $z \in \mathbb{C}^{+}$,
\begin{align}\label{Guassian-int}
\E&(S_{\mathbb{G}_{n'}}(z))-\E(S_{{\mathbb{G}}_{n',m}}(z))\notag\\
&= \frac{n'}{2n} \sum_{j=1}^{n} \sum_{k,\ell=n+1}^{n+N} \int_0^1
\Big ( {\mathbb{E}} (Z_{k-n}Z_{\ell-n}) - {\mathbb{E}}
({Z}_{k-n,m}{Z}_{\ell-n,m})\Big ) {\mathbb{E}} \Big( \frac{\partial^2}{\partial {x_{k,j}}\partial{x_{\ell,j}}}
f ( {\mathbf{g}} (t) ) \Big) dt \notag \\
& = \frac{N+n}{2n} \sum_{j=1}^{n} \sum_{k,\ell=1}^{N} \int_0^1
\Big ( {\mathbb{E}} (Z_{k}Z_{\ell}) - {\mathbb{E}}
({Z}_{k,m}{Z}_{\ell,m})\Big ) {\mathbb{E}} \Big( \frac{\partial^2}{\partial {x_{k+n,j}}\partial{x_{\ell+n,j}}}
f ( {\mathbf{g}} (t) ) \Big) dt
\end{align}
where, for $t \in[0,1]$, ${\mathbf{g}} (t)= \sqrt{N+n}\big(\sqrt{t} (\mathbb{G}_{n'})_{k,\ell} + \sqrt{1-t}
(\mathbb{G}_{n',m})_{k,\ell}\big)_{1 \leqslant \ell \leqslant k\leqslant n'}$ and $f$ is the function that allows us to write the Stieltjes transform of a symmetric matrix in terms of its entries (for a precise definition, see (49) in \cite{BaMePe}). Then, by using \eqref{covG} and \eqref{cov-Zk,m}, we write the following decomposition 
\begin{align}\label{deco-cavariances}
\E(Z_{k}Z_{\ell})-\E(Z_{k,m}Z_{\ell,m})
&= \E(X_{k}X_{\ell}) - \E(\bar{X}_{k,m}^{*}\bar{X}_{\ell,m}^{*})\notag \\
&= \E\big(X_{k}(X_{\ell} -\bar{X}_{\ell,m}^{*}) \big) 
+ \E\big( \bar{X}_{\ell,m}^{*}(X_k-\bar{X}_{k,m}^{*}  ) \big)
\end{align}
We shall decompose the left hand side of \eqref{Guassian-int} into two sums according to the decomposition \eqref{deco-cavariances} and treat them separately. By Lemma 13 in \cite{MePe}, we get for any two sequences $(a_{k})_{k}$ and $(b_{k})_{k}$ of real numbers
\begin{equation}\label{lemme-M.P.}
\sum_{j=1}^n \Big | \sum_{k,\ell=1}^N a_{k}b_{\ell}\frac{\partial^2}{\partial {x_{k+n,j}}\partial{x_{\ell+n,j}}}
f( {\mathbf{g}} (t) )\Big|
\leqslant \frac{C}{N+n}\Big (\sum_{k=1}^N a_{k}^{2}\sum_{\ell=1}^N b_{\ell}^{2}\Big )^{1/2}\,,
\end{equation}
where $C$ is a universal constant depending only on the imaginary part of $z$ and might change from a line to another. Applying \eqref{lemme-M.P.}  with $a_k=X_k$ and $b_{\ell}=X_{\ell} -\bar{X}_{\ell,m}^{*}$ and another time with $a_k=X_k- \bar{X}_{k,m}^{*}  $ and $b_{\ell}= \bar{X}_{\ell,m}^{*}$, we get for any  $z \in \mathbb{C}^{+}$,
\begin{equation}\label{bound-second-deriv-gaus-int}
\Big| \E(S_{\mathbb{G}_{n'}}(z))-\E(S_{\mathbb{G}_{n',m}}(z))\Big|^2
\leqslant \frac{C}{n^2}  \sum_{k=1}^{N} \big ( \|{X}_k\|^{2}_2 +   \|\bar{X}_{k,m}^{*}\|^2_2\big ) \sum_{\ell=1}^{N} \|X_{\ell}-\bar{X}_{\ell,m}^{*}\|^2_2 \,.
\end{equation}
In view of \eqref{Xk-Xk,mbar-norme2}, \eqref{X - M+}, \eqref{xk - Xk,m hat} and \eqref{norm Xbar - Xbar*}, we get 
\begin{align*}
\|X_{\ell}-\bar{X}_{\ell,m}^{*}\|^2 
&\leqslant
 2 \|X_{\ell}-\bar{X}_{\ell,m}\|^2 
+ 2\|\bar{X}_{\ell,m}-\bar{X}_{\ell,m}^{*}\|^2
\\ & \leqslant
16 \|X_0 -X_{0,m}\|_2^2 + 112 \| (|X_0|-M)_{+} \|_2^2 + 24 \| (|X_0|-L)_{+} \|_2^2 + 24 L^2 \beta_m \, ,
\end{align*}
where $L$ is a fixed positive real number strictly less than $M$ and not depending on $(n,m,M)$. Taking into account the stationarity of $(X_k)_k$ and \eqref{norm2 Xk,m}, we then infer  that 
\[
\Big| \E(S_{\mathbb{G}_{n'}}(z))-\E(S_{\mathbb{G}_{n',m}}(z))\Big|^2
\leqslant \frac{ C N^2}{n^2}  \Big(\|X_0 -X_{0,m}\|_2^2 +  \| (|X_0|-M)_{+} \|_2^2 +  \| (|X_0|-L)_{+} \|_2^2 +  L^2 \beta_m  \Big) \,,
\]
which converges to zero by letting first $n$, then $M$ followed by $m$ and finally $L$ tend to infinity. This ends the proof of the proposition.

\hfill$\square$

\subsection*{Acknowledgments}
I would like to thank the referee for carefully reading the manuscript and my thesis advisor Florence Merlev\`ede who was generous with her time, knowledge and assistance. I would also like to thank my co-advisor Emmanuel Rio for valuable discussions.

\end{document}